\documentclass[11pt,fleqn]{article}
\usepackage{amssymb}
\usepackage{amsfonts}
\usepackage{amsmath}
\usepackage{amsthm}
\usepackage{graphicx}
\usepackage{epsfig}
\usepackage{psfrag}
\usepackage{color}
\usepackage{dsfont}
\makeatletter
\xdef\@endgadget#1{{\unskip\nobreak\hfil\penalty50\hskip1em\hbox{}\nobreak
    \hfil#1\parfillskip=0pt\finalhyphendemerits=0\par}}
\def\@qedsymbol{${}_\blacksquare$}
\def\qed{\@endgadget{\@qedsymbol}}
\newtheorem{lemma}{Lemma}[section]

\newtheorem{example}[lemma]{Example}

\newtheorem{proposition}[lemma]{Proposition}
\newtheorem{remark}[lemma]{Remark}
\newcommand{\mR}{\mathbb{R}}

\DeclareMathOperator{\im}{im}

\DeclareMathOperator{\spa}{span}
\newcommand{\diag}{\mathrm{diag\,}}
\def\BibTeX{{\rm B\kern-.05em{\sc i\kern-.025em b}\kern-.08em
    T\kern-.1667em\lower.7ex\hbox{E}\kern-.125emX}}

\title{\LARGE \bf Modeling of Physical Network Systems\footnote{To appear in {\it Systems \& Control Letters}.}
}

\author{A.J. van der Schaft
\thanks{A.J. van der Schaft is with the Johann Bernoulli Institute for Mathematics and Computer
Science, and the Jan C. Willems Center for Systems and Control, University of Groningen, PO Box 407, 9700 AK, the
Netherlands
        {\tt\small A.J.van.der.Schaft@rug.nl}}
}
\date{}

\begin{document}

\maketitle
\thispagestyle{empty}
\pagestyle{empty}

\begin{abstract}
Conservation laws and balance equations for physical network systems typically can be described with the aid of the incidence matrix of a directed graph, and an associated symmetric Laplacian matrix. Some basic examples are discussed, and the extension to $k$-complexes is indicated. Physical distribution networks often involve a non-symmetric Laplacian matrix. It is shown how, in case the connected components of the graph are strongly connected, such systems can be converted into a form with balanced Laplacian matrix by constructive use of Kirchhoff's Matrix Tree theorem, giving rise to a port-Hamiltonian description. Application to the dual case of asymmetric consensus algorithms is given. Finally it is shown how the minimal storage function for physical network systems with controlled flows can be explicitly computed.
\end{abstract}

%\begin{keyword}
%%% keywords here, in the form: keyword \sep keyword
%physical network \sep Laplacian matrix \sep Matrix Tree theorem \sep port-Hamiltonian system \sep available storage.\\
%%% MSC codes here, in the form: 
%\MSC[2010]  93D20 \sep 94C15 \sep 05C50 \sep 15B51
%%% or \MSC[2008] code \sep code (2000 is the default)
%
%\end{keyword}

%\end{frontmatter}

% \linenumbers

\section{Introduction}
The topic of physical network systems has been always dear to Jan Willems' heart, from his early work on network synthesis and physical systems theory to his seminal work on dissipativity theory \cite{W1,W2}, and from the initial developments in behavioral theory to more recent 'educational' papers \cite{willems, willemsportsterminals}. I was often fortunate to witness these scientific developments from a close distance, and to be involved in penetrating discussions with Jan. Many of these animated debates centered around the 'right' and 'ultimate' definition of the basic concepts. Needless to say that my own ideas, including the ones presented in this paper, have been heavily influenced by Jan's.

The structure of the paper is as follows. In Section 3, after a recap of basic notions in algebraic graph theory in Section 2, I will discuss how {\it conservation laws} and {\it balance equations} for physical network are often naturally expressed in terms of the incidence matrix of a directed graph, and how this leads to a well-defined class of systems involving a {\it symmetric Laplacian matrix}. Next, in Section 4, attention will be directed to a more general class of physical network systems, of general distribution type, where the Laplacian matrix is {\it not} necessarily symmetric. Under the assumption of strong connectedness it will be shown how by means of Kirchhoff's Matrix Tree theorem the system can be constructively converted into a system with {\it balanced} Laplacian matrix, admitting a stability analysis similar to the symmetric case.
Section 5 is devoted to the analysis of {\it available storage} of passive physical network systems; a fundamental concept introduced in Jan Willems' seminal paper \cite{W1}. 
%In particular, it will be shown how explicit expressions relating directly to the graph structure can be obtained. 
Section 6 contains conclusions.

\section{Preliminaries about graphs}
We recall from e.g. \cite{bollobas, godsil} a few standard definitions and facts. 
A {\it graph} $\mathcal{G}(\mathcal{V},\mathcal{E})$, is defined by a set $\mathcal{V}$ of {\it vertices} (nodes) and a set $\mathcal{E}$ of {\it edges} (links, branches), where $\mathcal{E}$ is identified with a set of unordered pairs $\{i,j\}$ of vertices $i,j \in \mathcal{V}$. We allow for multiple edges between vertices, but not for self-loops $\{i,i \}$. By endowing the edges with an orientation we obtain a {\it directed graph}.
A directed graph with $n$ vertices and $m$ edges is specified by its $n \times m$ {\it incidence matrix}, denoted by $D$. Every column of $D$ corresponds to an edge of the graph, and contains exactly one $-1$ at the row corresponding to its tail vertex and one $+1$ at the row corresponding to its head vertex, while the other elements are $0$. In particular, $\mathds{1}^TD=0$ where $\mathds{1}$ is the vector of all ones. Furthermore, $\ker D^T = \spa \mathds{1}$ if and only if the graph is {\it connected} (any vertex can be reached from any other vertex by a sequence of, - undirected -, edges). In general, the dimension of $\ker D^T$ is equal to the number of connected components. A graph is {\it strongly connected} if any vertex can be reached from any other vertex by a sequence of directed edges.
For any diagonal positive semi-definite $m \times m$ matrix $R$ we define a {\it symmetric Laplacian matrix} of the graph as $L:=DRD^T$, where the positive diagonal elements $r_1, \cdots, r_m$ of the matrix $R$ are the weights of the edges. It is well-known \cite{bollobas} that $L$ is {\it independent} of the orientation of the graph.
%, and thus is associated with the undirected graph.

The vertex space \cite{hamgraphs} $\Lambda_0$ is defined as the set of all functions from the vertex set $\mathcal{V}$ to $\mR$. Obviously $\Lambda_0$ can be identified with $\mR^n$. The dual space of $\Lambda_0$ is denoted by $\Lambda^0$. Furthermore, the edge space $\Lambda_1$ is defined as the linear space of functions from the edge set $\mathcal{E}$ to $\mR$, with dual space denoted by $\Lambda^1$. Both spaces can be identified with $\mR^k$. It follows that the incidence matrix $D$ defines a linear map (denoted by the same symbol) $D: \Lambda_1 \to \Lambda_0$ with adjoint map $D^T: \Lambda^0 \to \Lambda^1$. 
Using these abstractions it is straightforward to extend the physical network dynamics described in this paper to other spatial domains than $\mR$. Indeed, for any linear space $\mathcal{R}$ (e.g., $\mathcal{R} = \mR^3$) we can define $\Lambda_0$ as the set of functions from $\mathcal{V}$ to $\mathcal{R}$, and $\Lambda_1$ as the set of functions from $\mathcal{E}$ to $\mathcal{R}$. In this case we can identify $\Lambda_0$ with the tensor product $\mR^n \otimes  \mathcal{R}$ and $\Lambda_1$ with the tensor product $\mR^k \otimes  \mathcal{R}$. Furthermore, the incidence matrix $D$ defines a linear map $D \otimes I : \Lambda_1 \to \Lambda_0$, where $I$ is the identity map on $\mathcal{R}$. In matrix notation $D \otimes I$ equals the Kronecker product of the incidence matrix $D$ and the identity matrix $I$. See \cite{hamgraphs} for further details.

\section{Physical network systems with symmetric Laplacian matrices}
%Additivity of energy function. Everything times n-dimensional, k edges. Use abstract notation for easy generalization.
%\newline
%\subsection{Conservation laws, balance equations, and symmetric Laplacians}
The structure of physical network dynamics is usually based on {\it conservation laws} and {\it balance equations}. Given a directed graph $\mathcal{G}$ with incidence matrix $D$ the basic way of expressing conservation laws is by equations of the form
\begin{equation}\label{cons}
Df +f_S=0,
\end{equation}
where $f \in \Lambda_1 \simeq \mathbb{R}^m$ is the vector of {\it flows} through the edges of the graph, and $f_S  \in \Lambda_0  \simeq \mathbb{R}^n$ is the vector of {\it injected flows} at the vertices. This can be regarded as a generalized form of Kirchhoff's current laws, in which case $f$ denotes the vector of currents through the edges of the electrical circuit graph, and $f_S$ are additional currents injected at the vertices of the circuit graph\footnote{Indeed, the presence of flows injected at the vertices is essential in Kirchhoff's original paper \cite{Kirchhoff}.}. Its restricted form is $Df=0$ (no flows/currents injected at the vertices). 

The injected flows at the vertices either correspond to external flows or to {\it storage} at the vertices, in which latter case there are state variables $x_i \in \mR$ (or, see above, $x_i$ belonging to a general linear space $\mathcal{R}$) associated to each $i$-th vertex, corresponding to $\dot{x} = -f_S$. This leads to the differential equations
\begin{equation}
\dot{x} = Df,
\end{equation}
expressing the basic {\it conservation laws} of the system: the sum of the incoming and outgoing flows through the edges incident to the $i$-th vertex is equal to the rate of storage at that vertex.

Often, the flows $f \in \Lambda_1$ through the edges are determined by {\it efforts} $e \in \Lambda^1$ associated to the edges, through a {\it resistive relation} of the form $f= - R(e)$, for some map $R: \Lambda^1 \to \Lambda_1$ satisfying $e^TR(e) \geq 0$, which is usually {\it diagonal} in the sense that its $j$-th component only depends on the effort $e_j$ associated to the $j$-th edge. The components of $e$ thus can be regarded as 'driving forces' for the flows $f$. In the linear case $f=-Re$ with $R$ a diagonal $n \times n$ matrix with nonnegative diagonal elements (in an electrical circuit context corresponding to conductances of resistors at the edges).

In many cases of interest, the effort variable $e_j$ corresponding to the $j$-th edge is an {\it across} variable which is determined by the {\it difference} of effort variables at the vertices incident to that edges, i.e.,
\begin{equation}
e = D^T e_S \, ,
\end{equation}
with $e_S \in \Lambda^0$ the vector of effort variables at the vertices. This corresponds to a {\it balance law} or an {\it equilibrium condition}: the driving force $e_j$ for the flow through the $j$-th edge is zero whenever the efforts at the vertices incident to this edge are equal. (In an electrical circuit $e_S$ corresponds to the voltage potentials at the vertices, and $e$ to the voltages across the edges.) 

Typically\footnote{The electrical circuit case is somewhat different, since resistors, inductors, and capacitors are all associated with the edges, and thus there is no storage at the vertices. Storage of charge at the vertices would correspond to {\it grounded capacitors}, with the ground node not included in the set of vertices.} the efforts $e_S$ at the vertices are determined by the state variables $x$ following
\begin{equation} 
e_S = \frac{\partial H}{\partial x}(x),
\end{equation} 
where $H: \Lambda_0 \to \mathbb{R}$ is the {\it total stored energy} at the vertices. Usually, $H$ is an {\it additive} energy function $H(x)= H_1(x_1) + \cdots + H_n(x_n)$. This leads to the equations
\begin{equation}\label{sym}
\dot{x} = -DRD^T \frac{\partial H}{\partial x}(x)
\end{equation}
The $n \times n$ matrix $L:= DRD^T$ is a {\it symmetric Laplacian} matrix, that is a symmetric matrix with nonnegative diagonal elements and nonpositive off-diagonal elements whose column and row sums are zero. Conversely, any symmetric Laplacian matrix can be represented as $L=DRD^T$ for some incidence matrix $D$ and positive diagonal matrix $R$. Clearly $L$ is positive semi-definite.

Equation (\ref{sym}) is the common form of physical network systems with energy storage confined to the vertices. (See \cite{hamgraphs} for other cases, in particular including energy storage associated to the edges.) They can be immediately seen to be in {\it port-Hamiltonian} form.
Recall, see e.g. \cite{vanderschaftmaschkearchive, NOW}, that port-Hamiltonian systems with inputs and outputs, in the absence of algebraic constraints and having linear energy-dissipating relations, are given by equations of the form
\begin{equation}
\begin{array}{rcl}
\dot{x} & = & \left[ \mathcal{J}(x) - \mathcal{R}(x) \right] \frac{\partial H}{\partial x}(x) + g(x)u \\[2mm]
y & = & g^T(x) \frac{\partial H}{\partial x}(x),
\end{array}
\end{equation}
with $\mathcal{J}(x)= - \mathcal{J}^T(x)$ and $\mathcal{R}(x)=\mathcal{R}^T(x) \geq 0$.
The system (\ref{sym}) is obviously port-Hamiltonian (without inputs and outputs) with  $\mathcal{J}(x) =0$ and $ \mathcal{R}(x)=L=DRD^T$.

\begin{example}[Mass-damper systems]
A paradigmatic example of the above scenario is a linear {\it mass-damper system}
\begin{equation}
\dot{p} = -DRD^T  M^{-1}p,
\end{equation}
with $p$ the vector of {\it momenta} of the masses associated to the vertices, $M$ the diagonal mass matrix, $R$ the diagonal matrix of damping coefficients of the dampers attached to the edges, and $H(p) = \frac{1}{2}p^T M^{-1}p$ the total kinetic energy of the masses. The vector of velocities $v=M^{-1}p$ converges to a vector in the kernel of $L=DRD^T$. In particular, if the graph is weakly connected the vector $v$ converges to a vector of the form $v^*\mathds{1}$, with $v^* \in \mathbb{R}$ (equal velocities).
For extensions to mass-spring-damper systems see \cite{hamgraphs}.
\end{example}

\begin{example}[Hydraulic networks]
Consider a hydraulic network between $n$ fluid reservoirs whose storage is described by the elements of a vector $x$. Mass balance corresponds to $\dot{x} = Df$ where $f \in \mathbb{R}^k$ is the flow through the $k$ pipes linking the reservoirs. Let each storage variable $x_i$ determine a pressure $\frac{\partial H_i}{\partial x_i}(x_i)$ for a certain energy function $H_i$. Assuming that the flow $f_j$ is proportional to the difference between the pressure of the head reservoir and the pressure of the tail reservoir this leads to the equations (\ref{sym}).
\end{example}

\begin{example}[Symmetric consensus algorithms] The equations (\ref{sym}) for $H(x) = \frac{1}{2} \|x\|^2$ reduce to $\dot{x} = -Lx, L=DRD^T$, which is the standard symmetric consensus protocol in continuous time, with weights given by the diagonal elements of $R$. In Section \ref{asymconsensus} we will pay attention to {\it asymmetric} consensus dynamics.
\end{example}

\subsection{Extension to higher-order complexes}
Directed graphs can be understood as the simplest case of $k$-complexes, namely as $1$-complexes where the incidence matrix $D$ is mapping edges to vertices (see also Chapter 13 of \cite{mesbahi}). A general $k$-complex is defined by a {\it sequence} of incidence (sometimes called 'boundary') operators
\[
\Lambda_k \overset{\partial_k}{\to} \Lambda_{k -1}
\overset{\partial_{k -1}}{\to}  \cdots \Lambda_1
\overset{\partial_1}{\to} \Lambda_0 
\]
with the property that
$
\partial_{j-1} \circ \partial_j = 0, \, j=2, \cdots , k.
$
The vector spaces $\Lambda_j, \,j=0,1 \cdots,k,$ are called
the spaces of $j$-chains. Each $\Lambda_j$ is generated by a
finite set of $j$-cells (like edges and vertices for graphs) in
the sense that $\Lambda_j$ is the set of functions from the
$j$-cells to $\mathbb{R}$. A typical example of a $k$-complex is the
triangularization of a $k$-dimensional manifold, with the
$j$-cells, $j=0,1, \cdots,k,$ being the sets of vertices, edges,
faces, etc..
Denoting the dual linear spaces by $\Lambda^j, \,j=0,1 \cdots,k,$
we obtain the following dual sequence
\[
\Lambda^0 \overset{d_1}{\to} \Lambda^1 \overset{d_2}{\to}
\Lambda^2  \cdots \Lambda^{k-1}  \overset{d_k}{\to} \Lambda^k
\]
where the adjoint maps $d_j \,j=0,1 \cdots,k,$ satisfy the analogous property
$
d_j \circ d_{j-1} =0, \, j=2, \cdots,k.
$
The elements of $\Lambda^j$ are called $j$-cochains.

Dynamics on the $k$-complex
can be defined in various ways, by defining resistive or energy-storing relations between the components of $\Lambda_j$
and $\Lambda^j, j=1, \cdots,k$; cf. \cite{vanderschaftmaschkeBosgra}. 

\begin{example}[Heat transfer on a $2$-complex] 
We will write the heat transfer in terms of the conservation of
internal energy. First we identify
the physical variables as chains and cochains of the given
$2$-complex. The components of the internal energy vector $u\in\Lambda^{2}$ denote the energy of each face. The heat conduction is given by the \emph{heat flux} $f \in\Lambda^{1}$ whose components equal the heat flux through every edge. Hence the basic conservation law (conservation of energy) is given as
\[
\frac{du}{dt}= d_{2}f 
\]
The thermodynamic properties are defined by Gibbs' relation, and
generated by the \emph{entropy function} $s: \Lambda^{2} \to \mathbb{R}$ as thermodynamic potential.
Since we consider transformations which are isochore and without
mass transfer, Gibbs' relation reduces to the definition of the
vector of intensive variables $e_u \in \Lambda_{2}$ which is (entropy-)\\conjugated to the
vector of extensive variables $u$:
\[
e_u =\frac{\partial s}{\partial u}(u)
\]
The components $e_u$ are equal to the reciprocal of the temperature at each $2$-face.

Since the temperature is varying over the faces, there is a \emph{thermodynamic driving force} vector 
$e \in\Lambda_{1}$
given by $e = \partial_{2}e_u$. 
By Fourier's law the heat flux is determined by the thermodynamic driving force vector as
\begin{equation}\label{eq:FourierDiscrete}
f = R(e_u)\, e,
\end{equation}
with $R(e_u)=R^T(e_u) \geq 0$ depending on the heat
conduction coefficients (note the opposite sign). 
The resulting system is a port-Hamiltonian system (with opposite sign of $R$ !), with vector of state
variables $x$ given by the internal energy vector $u$, and Hamiltonian
$s(u)$. It directly follows that the time-derivative of the entropy $s(u)$ satisfies
\[
\frac{ds}{dt} =  (\partial_2\frac{\partial
 s}{\partial u}(u))^T \,R(e_u)\, \partial_2\frac{\partial s}{\partial u}(u) = f^T\,R(e_u) \,f \geq 0
\]
expressing the fact that the entropy is monotonously {\it increasing}.

Note that the location of state variables corresponding to energy storage and resistive relations is somewhat orthogonal to (\ref{sym}), since the energy storing relation is between $\Lambda_2$ and $\Lambda^{2}$, and the resistive relation between $\Lambda_1$ and $\Lambda^{1}$, while in the case of (\ref{sym}) the energy storing relation is between $\Lambda_0$ and $\Lambda^{0}$, and the resistive relation between $\Lambda_1$ and $\Lambda^{1}$. This also leads to a different symmetric Laplacian matrix; see the discussion in \cite{vanderschaftmaschkeBosgra, seslija}.
\end{example}

\section{Physical network systems with non-symmetric Laplacian matrices}
On the other hand, not all physical network systems give rise to symmetric Laplacian matrices, at least not from the very start. In particular, it is not always the case that the effort variables $e$ associated to the edges (the driving forces for the flows $f$ through the edges) are given by the differences of effort variables associated to the incident vertices. For example, irreversible chemical reaction networks are not of this type; see e.g. \cite{rao, vdsJMC}. As a simpler example of the same type, let us consider linear {\it transportation networks}. Define the $m \times n$ transportation matrix $K$ as a matrix with $(j,i)$-th element equal to a positive weight constant $k_j$ if the transportation flow $f_j$ along the $j$-th edge is originating from vertex $i$ and equal to $f_j = k_jx_i$, and zero otherwise. Then the linear transportation network is represented by
\begin{equation}\label{lintransport}
\dot{x} = DKx,
\end{equation}
where the (generally asymmetric) matrix $L:=-DK$ has nonnegative diagonal elements and nonpositive off-diagonal elements, and furthermore satisfies $\mathds{1}^TL=0$. We will call such a matrix a {\it flow-Laplacian} matrix\footnote{In \cite{chapman} $L$ was called an {\it out-degree} Laplacian matrix.}. Conversely, any flow-Laplacian matrix $L$ can be represented as $L= -DK$ for some incidence matrix $D$ and a matrix $K$ consisting of nonnegative elements. An equivalent representation of a flow-Laplacian matrix is $L= \Delta - A$, where $A$ is the {\it adjacency matrix} of the directed graph with weights $k_j$ corresponding to the $j$-th edge, and where $\Delta$ is the diagonal matrix whose $i$-th element is equal to the sum of the weights of the edges originating from the $i$-th vertex.
\begin{remark}
A more general form of a linear transportation network is 
\begin{equation}
\dot{x} = Mx,
\end{equation}
where $M$ is an $n \times n$ {\it Metzler matrix}, that is, a matrix with {\it nonnegative off-diagonal} elements, see e.g. \cite{rantzer}. This still can be seen to define a system of the form (\ref{lintransport}), if we assume that $M$ is {\it diagonally dominant} in the sense that
\[
- m_{ii} \geq \sum_{i \neq j}m_{ij}, \quad i=1, \cdots, n,
\]
(implying in particular that $m_{ii} \leq 0, i=1, \cdots, n$). In this case an underlying graph structure can be defined as follows. The off-diagonal elements of $M$ define the off-diagonal elements of an adjacency matrix $A$ of a graph $\mathcal{G}$: if $m_{ji}$ is different from zero for $j \neq i$ then there is an edge from vertex $i$ to vertex $j$ with weight $m_{ji}$. We {\it augment} this graph $\mathcal{G}$ by a {\it sink vertex} and a corresponding additional row $(m_{n+1,1}, \cdots, m_{n+1,n})$ for its adjacency matrix $A$, with elements defined by
\[
m_{n+1,j} := - \sum_{i=1}^n m_{i,j}, \quad j=1, \cdots, n
\]
By adding the row vector $(m_{n+1,1}, \cdots, m_{n+1,n})$ as a last row to the matrix $M$, and finally adding an $(n+1)$-dimensional zero vector as last column, we obtain the negative of an $(n+1) \times (n+1)$ flow-Laplacian matrix $L$. (Note that the sink vertex of the augmented graph takes care of the {\it surplus} of flow in the original graph $\mathcal{G}$.)
\end{remark}

\subsection{Dynamics of distribution networks}
As motivated above, a large class of distribution networks (including subclasses of transportation networks) is of the form
\begin{equation}\label{nonsym}
\dot{x} = - L\frac{\partial H}{\partial x}(x),
\end{equation}
where $L$ is a flow-Laplacian matrix, and $H(x)= H_1(x_1) + \cdots + H_n(x_n)$ is an additive energy function. 

What is the structure of the dynamics (\ref{nonsym}), and how can it be analyzed ? Note that a flow-Laplacian matrix $L$ does not necessarily satisfy $L + L^T \geq 0$, and consequently the energy $H$, even if it is bounded from below, is not necessarily a Lyapunov function for (\ref{nonsym}). In fact\footnote{See e.g. \cite{cortes}, where the proof is given for a matrix $L$ such that $L^T$ is a flow-Laplacian matrix.}
\begin{proposition}
The flow-Laplacian matrix $L$ satisfies $L + L^T \geq 0$ if and only if it is {\it balanced}; that is, not only $\mathds{1}^TL=0$ (column sums zero) but also $L\mathds{1}=0$ (row sums zero). 
\end{proposition}
The main aim of this subsection is to show how a flow-Laplacian matrix can be {\it transformed} into a balanced one, provided the connected components of the graph are all {\it strongly connected}. (Closely related developments can be found in \cite{lewis}.) Furthermore, we show how to do this in a constructive way by employing a general form of Kirchhoff's Matrix Tree theorem\footnote{This theorem goes back to the classical work of Kirchhoff on resistive electrical circuits \cite{Kirchhoff}; see \cite{bollobas} for a succinct treatment, and \cite{gunawardena,  vdsJMC} for an account in the context of chemical reaction networks.}, also providing additional insights. 

First assume that the graph under consideration is connected, implying that $\dim \ker L =1$. Kirchhoff's Matrix Tree theorem tells us how we can find a nonnegative vector $\sigma$ in the kernel of $L$ as follows.
Denote the $(i,j)$-th cofactor of $L$ by $C_{ij}=(-1)^{i+j}M_{i,j}$, where $M_{i,j}$ is the determinant of the $(i,j)$-th minor of $L$, which is the matrix obtained from $L$ by deleting its $i$-th row and $j$-th column. Define the adjoint matrix $\mathrm{adj}(L)$ as the matrix with $(i,j)$-th element given by $C_{ji}$. It is well-known that
\begin{equation}\label{adjoint}
L \cdot \mathrm{adj}(L) = (\det{L})I_n =0
\end{equation}
Furthermore, since $\mathds{1}^TL=0$ the sum of the rows of $L$ is zero, and hence by the properties of the determinant function $C_{ij}$ does not depend on $i$; implying that $C_{ij} = \sigma_j, \, j=1, \cdots, n$. Hence by defining $\sigma := (\sigma_1, \cdots, \sigma_n)^T$, it follows from (\ref{adjoint}) that $L\sigma=0$. Furthermore, cf. \cite[Theorem 14 on p.58]{bollobas}, $\sigma_i$ is equal to {\it the sum of the products of weights of all the spanning trees of} $\mathcal{G}$ {\it directed towards} vertex $i$. In particular, it follows that $\sigma_j \geq 0, j=1, \cdots,c$. In fact, $\sigma \neq 0$ if and only if $\mathcal{G}$ has a spanning tree. Furthermore, since for every vertex $i$ there exists at least one spanning tree directed towards $i$ if and only if the graph is strongly connected, we conclude that $\sigma \in \mathbb{R}^n_+$ if and only if the graph is {\it strongly connected}. 
\begin{remark}
The {\it existence} (not the explicit {\it construction}) of $\sigma \in \mathbb{R}^n_+$ satisfying $L \sigma=0$ already follows from the Perron-Frobenius theorem; exploiting the fact that the off-diagonal elements of $-L:=DK$ are all nonnegative \cite[Lemma V.2]{sontag}.
%\footnote{This implies that there exists a real number $\alpha$ such that $-L + \alpha I_m$ is a matrix with all elements nonnegative. Since the set of eigenvectors of $-L$ and $-L + \alpha I_m$ are the same, and moreover by $\mathds{1}^TL=0$ there cannot exist a positive eigenvector of $-L$ corresponding to a non-zero eigenvalue, the application of Perron-Frobenius to $-L + \alpha I_m$ yields the result; see \cite[Lemma V.2]{sontag} for details.}.
\end{remark}
%\begin{example}\label{excyclic}
%Consider the cyclic reaction network
%\begin{tabular}{c c c}
%& $C_3$ & \\
%& {\rotatebox[origin=c]{45}{$\xleftrightharpoons[k_3^-]{\ k_3^+ \ }$}} \ {\rotatebox[origin=c]{-45}{$\xleftrightharpoons[k_2^-]{\ k_2^+ \ }$}} & \\
%& $C_1$ \ \ \ \ \ \ $\xrightleftharpoons[k_1^-]{\ k_1^+ \ }$  \ \ \ \ \ \ $C_2$ &
%\end{tabular}
%\vspace{0.3cm}
%
%\noindent in three (unspecified) complexes $C_1, C_2, C_3$. The Laplacian matrix is given as
%\[
%L = \begin{bmatrix} k_1^+ + k_3^- & - k_1^- & - k_3^+ \\
%- k_1^+ & k_1^- + k_2^+ & -k_2^- \\
%-k_3^- & - k_2^+ & k_3^+ + k_2^-
%\end{bmatrix}
%\]
%By Kirchhoff's Matrix Tree theorem the corresponding vector $\rho$ satisfying $L \rho=0$ is given as
%\[
%\rho = \begin{bmatrix} k_2^+k_3^+ + k_1^-k_3^+ + k_1^-k_2^- \\
%k_1^+k_3^+ + k_1^+k_2^- + k_2^-k_3^- \\
%k_1^+k_2^+ + k_2^+k_3^- + k_1^-k_3^- 
%\end{bmatrix},
%\]
%where each term corresponds to one of the three weighted spanning trees pointed towards the three vertices.
%
%\end{example}
In case the graph $\mathcal{G}$ is not connected the same analysis can be performed on each of its connected components. Hence if all connected components of $\mathcal{G}$ are strongly connected, Kirchhoff's Matrix Tree theorem provides us with a constructive way to obtain a vector $\sigma \in \mathbb{R}^n_+$ such that $L \sigma = 0$. It immediately follows that the transformed matrix
\begin{equation}
\mathcal{L} := L \Sigma,
\end{equation}
where $\Sigma$ is the $n \times n$-dimensional diagonal matrix with diagonal elements $\sigma_1, \cdots, \sigma_n,$ is a balanced Laplacian matrix\footnote{Conversely, it is known that a graph with balanced Laplacian matrix is connected if and only if it is strongly connected.}. Hence we can rewrite the dynamics (\ref{nonsym}) as
\begin{equation}
\dot{x} = - \mathcal{L} \frac{\partial \mathcal{H}}{\partial x}(x)
\end{equation}
with the additive transformed energy function $\mathcal{H}(x) = \mathcal{H}_1(x_1) + \cdots + \mathcal{H}_n(x_n)$ defined by $\sigma_i \mathcal{H}_i(x_i) = H_i(x_i), \, i=1, \cdots, n.$ By splitting $\mathcal{L}$ into its skew-symmetric and symmetric part $\mathcal{L}= - \mathcal{J} + \mathcal{R}$ with $\mathcal{R}=\mathcal{R}^T \geq 0$, we obtain the port-Hamiltonian representation $\dot{x} =  [ \mathcal{J} - \mathcal{R} ]\frac{\partial \mathcal{H}}{\partial x}(x)$. In case the original energy function $H$ is bounded from below also $\mathcal{H}$ is bounded from below, and thus defines a Lyapunov function, since 
\[
\frac{d}{dt} \mathcal{H}(x) = -\frac{\partial^T \mathcal{H}}{\partial x}(x) \mathcal{L} \frac{\partial \mathcal{H}}{\partial x}(x) \leq 0
\]

\subsection{Interconnection of passive systems}
A direct extension of the previous theory concerns interconnection of passive systems. Consider $N$ nonlinear systems of the form
\begin{equation}
\begin{array}{rcl}
\dot{x}_i & = & f_i(x_i) + g_i(x_i)u_i , \quad i=1, \cdots,N , \\[2mm]
y_i & = & h_i(x_i)
\end{array}
\end{equation}
with scalar inputs and outputs,
which are assumed to be {\it passive}; that is, there exist (differentiable) storage function $S_i(x_i) \geq 0$ satisfying
\begin{equation}
\frac{\partial S_i}{\partial x_i}(x_i) \leq h_i(x_i)u_i, \quad i=1, \cdots,N
\end{equation}
These systems are associated to the $N$ vertices of a directed graph $\mathcal{G}$, which are assumed to be linked via the $M$ edges of the graph by
\begin{equation}
u = -Ly + v,
\end{equation}
where $u$ and $y$ are the vectors with subvectors $u_i, y_i, i=1, \cdots,M$, and $L$ is a flow-Laplacian matrix of the graph. Defining $S:= S_1 + \cdots + S_N$ the interconnected system obviously satisfies
\begin{equation}
\frac{d}{dt}S \leq y^Tu = - y^TLy + y^Tv
\end{equation}
In case $L$ is {\it balanced} it follows that $y^TLy \geq 0$, showing that the interconnected system is passive with respect to the inputs $u_1,\cdots,u_M$ and outputs $y_1,\cdots,y_M$ with storage function $S$. This observation was one of the starting points of \cite{chopra} and subsequent papers.

However, if $L$ is not balanced, then in general the term $y^TLy$ is {\it not} $\geq 0$, and hence $S$ does not define a storage function for the interconnected system. On the other hand, assume that the connected components of $\mathcal{G}$ are {\it strongly connected}, then by the application of Kirchhoff's Matrix Tree theorem as above there exists a vector $\sigma = (\sigma_1, \cdots,\sigma_N)^T$ such that $L \sigma = 0$. By defining the corresponding {\it weighted} combination
\begin{equation}
S^{\sigma}(x_1, \cdots, x_N):= \frac{1}{\sigma_1} S_1(x_1) + \cdots + \frac{1}{\sigma_N} S_N(x_N)
\end{equation}
it follows that
\begin{equation}
\begin{array}{rcl}
\frac{d}{dt}S^{\sigma} & \leq & \frac{1}{\sigma_1} y_1 u_1 + \cdots + \frac{1}{\sigma_N} y_N u_N = y^T \Sigma^{-1} u \\[2mm]
&= &- y^T\Sigma^{-1} L y + y^T \Sigma^{-1}v \\[2mm]
& = & - y^T\Sigma^{-1} \mathcal{L} \Sigma^{-1}y + y^T \Sigma^{-1}v
\end{array}
\end{equation}
where, as above, $\Sigma = \diag (\sigma_1, \cdots, \sigma_N)$. 
Since the matrix $\mathcal{L} = L \Sigma$ is {\it balanced} $y^T\Sigma^{-1} \mathcal{L} \Sigma^{-1}y\geq 0$, and hence the interconnected system is passive with respect to the new inputs $v_1, \cdots, v_M$ and scaled outputs $ \frac{1}{\sigma_1}y_1, \cdots,  \frac{1}{\sigma_N}y_N$.
% (or scaled inputs $ \frac{1}{\sigma_1}v_1, \cdots,  \frac{1}{\sigma_N}v_N$ and original outputs $y_1, \cdots, y_N)$.

\subsection{Digression on asymmetric consensus algorithms}\label{asymconsensus}
The situation considered before is {\it dual} to non-symmetric {\it consensus algorithms} in continuous time; see e.g. \cite{mesbahi}. In this case one considers a multi-agent system with imposed dynamics
\begin{equation}
\dot{x}_i = - \sum_{j \in \mathcal{N}_i}a_{ij} (x_i - x_j), \quad i=1, \cdots, n,
\end{equation}
for certain constants $a_{ij} \geq 0$, where $\mathcal{N}_i \subset \{1, \cdots,n\}$ consists of all vertices from which there is an edge directed towards vertex $i$. In terms of the total vector $x=(x_1, \cdots,x_n)^T$ this can be written as
\begin{equation}\label{consensus}
\dot{x} = - L_cx,
\end{equation}
where $L_c$ is the $n \times n$ matrix with $(i,j)$-th off-diagonal element $-a_{ij}$, and with $i$-th diagonal element $\sum_j a_{ij}, i=1, \cdots,n$. Hence the matrix $L_c$ has nonnegative diagonal elements and nonpositive off-diagonal elements, and has zero {\it row sums}, that is, $L_c\mathds{1}=0$. Equivalently, $L_c^T$ is a flow-Laplacian matrix as considered before.

Application of Kirchhoff's Matrix Tree theorem to $L_c$ now yields the existence of a row vector $\sigma = (\sigma_1, \cdots, \sigma_n)$ such that $\sigma L_c=0$, where $\sigma_i$ is the {\it sum of the products of weights along spanning trees directed from} vertex $i$. What are the implications of this ? It means that along the consensus dynamics (\ref{consensus})
\begin{equation}
\frac{d}{dt} \sum_j \sigma_jx_j =0,
\end{equation}
thus defining a {\it conserved quantity}. This conserved quantity is non-trivial if there exists at least one spanning tree, which is \cite{mesbahi} the necessary and sufficient condition for convergence of the dynamics (\ref{consensus}) to consensus, that is, to a vector $c \mathds{1}$, where $c$ is the consensus value corresponding to the initial condition $x(0)$. In this case, reorder the vertices in such a way that for each of the first $k$ vertices there exists a spanning tree directed from that vertex, and for none of the last $n-k$ vertices. This means that $\sigma_j >0, j=1, \cdots,k,$ and $\sigma_j=0, j=k+1, \cdots,n$. It follows that the consensus value $c$ (depending on $x(0)=(x_1^0, \cdots,x_n^0)^T$) is determined by
\begin{equation}
c= \frac{\sigma_1x_1^0 + \cdots + \sigma_kx_k^0}{\sigma_1 + \cdots + \sigma_k}
\end{equation}
In particular, the consensus value is independent of the values of the initial condition of the state variables of those vertices from which there is no spanning tree directed from this vertex (and thus no global flow of information stemming from this vertex).

In case the graph is {\it strongly connected} it means that all elements of $\sigma$ are positive, and thus we can define, as in the case of a flow-Laplacian matrix, the positive diagonal matrix $\Sigma$ with diagonal elements $\sigma_1, \cdots, \sigma_n$, and rewrite the consensus dynamics (\ref{consensus}) as
\begin{equation}
\Sigma \dot{x} = - \Sigma L_c x = - \mathcal{L}_cx,
\end{equation}
where the Laplacian matrix $\mathcal{L}_c := \Sigma L_c$ is balanced. It follows that along this dynamics
\[
\frac{d}{dt}\frac{1}{2}x^T\Sigma x = - x^T \mathcal{L}_cx \leq 0
\]
showing convergence to consensus. Hence, the consensus dynamics (\ref{consensus}) admits an additive Lyapunov function. It also implies that asymmetric consensus dynamics can be written into port-Hamiltonian form
\[
\dot{z} = [ \mathcal{J} - \mathcal{R} ] \, \Sigma^{-1} z
\]
with $z:= \Sigma \dot{x}$, Hamiltonian $H(z)= \frac{1}{2}z^T \Sigma^{-1}z$, and $\mathcal{L}_c = - \mathcal{J} + \mathcal{R}$ the decomposition into skew-symmetric and positive semi-definite symmetric part.

\begin{remark}
Note that the 'consensus Laplacian matrix' $L_c$ can be alternatively written as $L_c=JD^T$ for a certain matrix $J$ containing the coefficients $l_{ij}$. It follows that from a coordinate-free point of view $L_c: \Lambda^0 \to \Lambda^0$, while $L: \Lambda_0 \to \Lambda_0$ for a flow-Laplacian matrix $L$. Hence, the state vector $x$ for physical networks is an element of $\Lambda_0$, the vertex space, while the state vector $x$ for consensus dynamics is an element of the dual vertex space $\Lambda^0$ (a linear function $z \mapsto x^Tz$ for $z \in \Lambda_0$, or, alternatively, a vector of potentials).
\end{remark}

\section{Available storage of passive physical network systems}
%How to extend this to mass-spring systems, or mass-damper systems ?\\
%Can we have systems on different graphs having the same input-output behavior ?\\
%What is damping going to change ?\\
%When is uniqueness of the internal storage guaranteed ?\\
%\\
%
Consider the following system on a directed graph $\mathcal{G}$
\begin{equation}\label{simplegraph}
\begin{aligned}
\dot{x} & = Du, \quad \quad u \in \mathbb{R}^m, \ x \in \mathbb{R}^n, \\
y & = D^T \frac{\partial H}{\partial x}(x), \quad y \in \mathbb{R}^m,
\end{aligned}
\end{equation}
where $D$ is the incidence matrix of $\mathcal{G}$, and $H(x)$ is the Hamiltonian function (physical energy). This corresponds to a system with 'flow source' at every edge, and outputs which are the differences of the potentials $\frac{\partial H}{\partial x_i}(x)$ at each $i$-th vertex. Throughout this section we will assume\footnote{Without loss of generality, since otherwise the analysis can be repeated for every connected component of the graph.} that the graph is {\it connected}, or equivalently $\ker D^T = \spa \mathds{1}$. 

If we assume $H$ to be non-negative then $H$ defines a storage function, and the system is passive. On the other hand, as we will show in this section\footnote{For an earlier version of the first part of this section see \cite{NOW}.}, the {\it minimal} storage function for the system will be always strictly smaller than $H$.
Recall from the seminal papers \cite{W1, W2} of Willems on dissipativity theory that the minimal storage function is equal to the available storage $S_a$ given as
\begin{equation}
S_a(x) = \sup - \int_0^{\tau} y^T(t)u(t) \, dt,
\end{equation}
where we consider the supremum over all $\tau \geq 0$ and all input functions $u: [0,\tau] \to \mathbb{R}^m$, and where $y: [0,\tau] \to \mathbb{R}^m$ is the output resulting from the input function $u: [0,\tau] \to \mathbb{R}^m$ and initial condition $x(0) = x$. 

Consider first the case $H(x)= \frac{1}{2} \|x \|^2$. Noting that
\begin{equation}\label{ava}
\begin{array}{l}
\int_0^{\tau} y^T(t)u(t) \, dt  = \int_0^{\tau} x^T(t)Du(t) \, dt =\\[2mm] 
 \int_0^{\tau} x^T(t) \dot{x}(t) \, dt = \frac{1}{2} \|x(\tau) \|^2 - \frac{1}{2} \|x(0) \|^2,
\end{array}
\end{equation}
we see that the available storage in this case is given as
\begin{equation}\label{availablegraph}
S_a(x) = \sup_{x(\tau)} \left( \frac{1}{2} \|x \|^2 - \frac{1}{2} \|x(\tau) \|^2 \right),
\end{equation}
where we take the supremum over all $\tau \geq 0$ and all possible states $x(\tau)$ resulting from input functions $u: [0,\tau] \to \mathbb{R}^m$. By connectedness of the graph, we know that from $x(0)=x$ we can reach, by choosing the input function suitably, any state $x(\tau)$ satisfying
\begin{equation}\label{conserved}
\mathds{1}^T x(\tau) = \mathds{1}^T x.
\end{equation}
Hence the available storage $S_a(x)$ is given by (\ref{availablegraph}) where the supremum is taken over all states $x(\tau)$ satisfying (\ref{conserved}). This corresponds to minimizing $\frac{1}{2} \|x(\tau) \|^2$ over all $x(\tau)$ satisfying (\ref{conserved}), having the solution
\begin{equation}
x(\tau) = \frac{1}{n} \mathds{1}^T x \mathds{1},
\end{equation}
Thus the available storage $S_a$ is given by the explicit expression
\begin{equation}
S_a(x) = \frac{1}{2} \|x \|^2 - \frac{1}{2} \left(\frac{1}{n} \mathds{1}^T x\right)^2 \| \mathds{1}\|^2 = \frac{1}{2} x^T \left(I_n - \frac{1}{n} \mathds{1}\mathds{1}^T\right)x
\end{equation}
We conclude that for all initial conditions $x(0)=x$ which are such that $\mathds{1}^Tx \neq 0$ the available storage $S_a(x)$ is typically smaller than the Hamiltonian $\frac{1}{2} \|x \|^2$. The reason is that, since the system $\dot{x} = Du$ is not controllable, it is not possible to drive every initial state to the origin; the position where the energy $H(x) = \frac{1}{2} \|x \|^2$ is zero. Instead, by extracting the maximal amount of energy the system is brought from state $x$ to a state $x^*$ with $x^*_1= \cdots = x^*_n$, satisfying $x^*_1 + \cdots + x^*_n = x_1 + \cdots + x_n$. Note that as a consequence, $H(x)=\frac{1}{2} \|x \|^2$ in (\ref{simplegraph}) may be replaced by $S_a$.

Also note that the matrix $I_n - \frac{1}{n} \mathds{1}\mathds{1}^T$ defines a symmetric weighted Laplacian matrix for an extended graph; namely the {\it complete graph} for the vertices of the original graph\footnote{A graph is {\it complete} if there is an edge between every pair of vertices.}.

The above analysis can be extended to any system (\ref{simplegraph}) for which the Hamiltonian $H$ is strictly convex and bounded from below. Indeed, for such an $H$ the available storage can be seen to be
\begin{equation}\label{storage1}
S_a(x) = H(x) - H(x^*(x)),
\end{equation}
where $H(x^*(x))$ is the solution of minimizing $H(x^*)$ over all $x^* \in \mathbb{R}^n$ satisfying $\mathds{1}^Tx^* = \mathds{1}^Tx$. Equivalently, this amounts to the minimization of $H(x^*) + \lambda (\mathds{1}^Tx^* - \mathds{1}^Tx)$ over $x^*$ and the Lagrange multiplier $\lambda \in \mathbb{R}$, yielding the minimizer $x^*(x)$ as the solution of the equations
\begin{equation}\label{opt}
\begin{array}{c}\displaystyle
\frac{\partial H}{\partial x_1}(x^*(x)) = \cdots = \frac{\partial H}{\partial x_n}(x^*(x)), \\[2mm]
x_1^*(x) +  \cdots + x_n^*(x) = x_1 + \cdots + x_n.
\end{array}
\end{equation}
As in the case $H(x) = \frac{1}{2} \|x\|^2$, the expression for the available storage is independent of the graph (as long as it is connected). Note furthermore that the first line of (\ref{opt}) can be interpreted as a {\it consensus} condition on $\frac{\partial H}{\partial x_1}, \cdots,  \frac{\partial H}{\partial x_n}$. 
\begin{example}
Consider a system of $n$ point masses $m_1, \cdots, m_n$ in $\mathbb{R}$ with state variables being the momenta $p_1, \cdots, p_n$, and with Hamiltonian equal to the {\it kinetic energy}
\[
H(p) =  \sum_{i=1}^n \frac{p_i^2}{2m_i}
\]
The available storage can be computed
\[
S_a(p) = \frac{1}{2}\sum_{i <j} \frac{m_im_j}{m_1 + \cdots + m_n} \left( \frac{p_i}{m_i} - \frac{p_j}{m_j} \right)^2.
\]
This quantity was called the {\it motion energy} in \cite{willemscdc13}. It amounts to the maximal energy which can be extracted from the system by applying forces $F_1, \cdots, F_n$ satisfying $\sum_{j=1}^nF_j =0$, or equivalently (since $\mathds{1}^TD=0$)
\[
\dot{p} = F= Du,
\]
where $F$ is the vector with components $F_1, \cdots, F_n$, and $D$ is the incidence matrix of the complete graph with vertices corresponding to the masses\break $m_1, \cdots, m_n$. Note that the same available storage results for any incidence matrix $D$ corresponding to a connected graph. As a result of extracting the maximal energy, the system will end up in a consensus state $v_1 = \cdots = v_n$, with $v_i = \frac{p_i}{m_i}$ the velocities of the point masses.

This can be readily extended to point masses in $3$-dimensional Euclidean space with $p_i \in \mathbb{R}^3$; replacing the expression 
\begin{equation*}
\left( \frac{p_i}{m_i} - \frac{p_j}{m_j} \right)^2\ \text{with} \ \ \left\| \frac{p_i}{m_i} - \frac{p_j}{m_j} \right\|^2. 
\end{equation*}
\end{example}
As a general remark we mention that contrary to the Hamiltonian function, the available storage is {\it not} necessarily the sum of the available storages of the individual subsystems, as was already noted in \cite{willemscdc13}. 
A simple example is provided by the juxtaposition of two systems each consisting of two masses. The sum of the energies which can be extracted from the two systems separately (by applying for each system two external forces whose sum is zero) is strictly smaller than the amount of energy which can be extracted from the four masses by applying four external forces whose sum is zero.

\subsection{Generalization}
How to generalize the above results to passive physical network systems, more general than (\ref{simplegraph}) ? An interesting case is where the edges are split into a subset of 'flow sources' and a complementary subset of edges corresponding to a resistive relation. This corresponds to splitting the incidence matrix $D$ as $D = [D_s  \, D_u]$, where $D_s$ corresponds to the flow sources and $D_u$ to the remaining edges. For linear resistive relations this yields the system description
\begin{equation}\label{generalized}
\begin{array}{rcl}
\dot{x} & = & -D_uRD_u^T \frac{\partial H}{\partial x}(x) + D_s u \\[2mm]
y & = & D_s^T\frac{\partial H}{\partial x}(x)
\end{array}
\end{equation}
for a certain positive diagonal matrix $R$. This can be seen to be a generalization both of (\ref{simplegraph}) (by adding resistive relations) and of (\ref{sym}) (by adding flow sources).

Throughout we will assume that the system (\ref{generalized}) is {\it controllable} restricted to each affine space $x(0) + \im D$. For the simplest case $H(x) = \frac{1}{2} \|x\|^2$ this can be checked using the Kalman controllability condition on the pair $(D_uRD_u^T,D_s)$: the smallest subspace containing $\im D_s$ and invariant under the partial Laplacian matrix $D_uRD_u^T$ should be equal to $\im D$. It can be verified that for generic values of the diagonal elements of $R$ this weak form of controllability is guaranteed.
%\footnote{Of course, specific physical examples may still be outside this generic class.}.

In case of a general Hamiltonian $H$ the expression (\ref{ava}) for the delivered power generalizes to
\begin{equation}\label{ava1}
\begin{array}{l}
\int_0^{\tau} y^T(t)u(t) dt  = \int_0^{\tau} \frac{\partial^T H}{\partial x}(x(t))D_su(t) \, dt =\\[2mm]
 \int_0^{\tau} \frac{\partial^T H}{\partial x}(x(t)) \, [\dot{x}(t) +  D_uRD_u^T \frac{\partial H}{\partial x}(x(t))] \,dt = \\[2mm] 
 H(x(\tau)) - H(x(0)) + \int_0^{\tau}\frac{\partial^T H}{\partial x}(x(t)) D_uRD_u^T \frac{\partial H}{\partial x}(x(t)) \, dt
\end{array}
\end{equation}
Hence the available storage $S_a(x)$ is bounded from above by the same expression $H(x) - H(x^*(x))$ as obtained in (\ref{storage1}), with $x^*(x)$ satisfying (\ref{opt}). Furthermore, whenever $x^*(x)$ is satisfying (\ref{opt}), the dissipated power\break $\frac{\partial^T H}{\partial x}(x^*x) D_uRD_u^T \frac{\partial H}{\partial x}(x^*(x))$ is equal to zero. Hence, starting from any initial state $x(0)=x$, the system may be steered to a state $x(\tau)=x^*(x)$ satisfying (\ref{opt}) and $\mathds{1}^Tx^*(x) = \mathds{1}^Tx(0)$ in such a way that the dissipated energy $\int_0^{\tau}\frac{\partial^T H}{\partial x}(x(t)) D_uRD_u^T \frac{\partial H}{\partial x}(x(t)) \, dt$ is arbitrarily small. Therefore, as for the system (\ref{simplegraph}), the available storage $S_a(x)$ is actually {\it equal} to (\ref{storage1}) with $x^*(x)$ satisfying (\ref{opt}).

\section{Conclusions}
Physical network systems are motivated by a variety of application areas, from circuit theory, mechanical networks, hydraulic systems and power networks, to systems biology. Furthermore, they naturally arise from finite-dimensional modelling of systems of a distributed-parameter nature\break \cite{vdsmaschkeJGM}, either by lumped systems modelling or by structure-preserving spatial discretization; see e.g. \cite{seslija1, seslija}. 

A major recent impetus to the study of complex physical network systems is the systematic use of concepts and tools from algebraic graph theory. Remarkably, this also constitutes a return to the origin of the subject: Kirchhoff's classical paper on resistive circuits \cite{Kirchhoff} can be regarded at the same time as one of the starting points of algebraic graph theory.

The blending of physics, systems theory, algebra and geometry, continues to provide inspiration; very much in the spirit of the work of Jan Willems.

\medskip



%% References


\bibliographystyle{elsarticle-num}



\end{document}